\documentclass{comnet}

\usepackage{graphicx,amsmath,amsthm,amssymb,latexsym,amsfonts,color,
dsfont}
\usepackage{subfig}

\begin{document}

\title{Structural analysis of water networks}

\shorttitle{Structural analysis of water networks} 
\shortauthorlist{M.~Benzi, I.~Daidone, C.~Faccio, L.~Zanetti-Polzi} 

\author{
\name{Michele Benzi}
\address{Scuola Normale Superiore, Piazza dei Cavalieri, 7, 56126, Pisa, Italy}
\name{Isabella Daidone}
\address{Department of Physical and Chemical Sciences, University of L'Aquila, via Vetoio
(Coppito 1), L'Aquila, 67010, Italy	}
\name{Chiara Faccio$^*$}
\address{Scuola Normale Superiore, Piazza dei Cavalieri, 7, 56126, Pisa, Italy\email{$^*$Corresponding author: chiara.faccio@sns.it}}
\and
\name{Laura Zanetti-Polzi}
\address{Center S3, CNR-Institute of Nanoscience, Via Campi 213/A, Modena, 41125, , Italy}}

\maketitle

\begin{abstract}
{Liquid water, besides being fundamental for life on Earth, has long fascinated scientists due to several anomalies.
Different hypotheses have been put forward to explain these peculiarities. The most accredited one foresees the presence in the supercooled region of two phases at different densities: the low-density liquid phase and the high-density liquid phase. In our previous work [Faccio et al., J. Mol. Liq. 355 (2022): 118922], we showed that it is possible to identify these two forms 
in water networks through a computational approach based on molecular dynamics simulation and on the calculation of the total communicability of the associated graph, 
in which the nodes correspond to water molecules and the edges represent the connections (interactions) between molecules. In this paper, we present a more 
in-depth investigation of the application of graph-theory based approaches to the analysis of the structure of water networks. In particular, we investigate different connectivity and centrality measures and we report on the use of a variety
of global metrics aimed at giving a topological and geometrical characterization of liquid water.}
{liquid water, LDL/HDL phases, graph structure, centrality measures, global metrics.}
\end{abstract}

\section{Introduction}\label{sec1}

The mathematical theory of graphs has long found important applications in chemistry.  In its most basic form it has been used, for example,  to construct mathematical models of molecules and their interactions. A host of graph-theoretic concepts and results can be applied
to analyze chemical compounds and understand their properties, and thousands of papers have been published in this field; see, e.g., \cite{Balaban}, or \cite{Yirik_et_al}
for a more recent survey.

Only recently, concepts and techniques from the field of network science have been introduced to chemical physics with the goal of gaining a deeper understanding of collective phenomena involving molecular networks. This is in contrast with other fields, like biology, where network science tools have been in heavy use for a number of years, for example in the study of protein and gene networks, see  \cite{gursoy2008topological,barabasi2004network}.  In the present paper we apply concepts from network science, together with molecular dynamics (MD) tools, to advance our understanding of the structure of liquid water through the analysis of networks of water molecules under different conditions. The paper can be seen as a follow-up to our previous work \cite{Benzi4}, but it is self-contained and can be read independently of it.

In this paper we are interested, primarily, in the 
problem of characterizing a liquid phase in a molecular simulation. In particular, we focus here on the differentiation between two forms of liquid water differing in density: the low-density liquid (LDL) and high-density liquid (HDL). The existence of a first-order liquid–liquid phase transition (LLPT) between these
two forms in the deeply supercooled region at high pressure has been hypothesized on the basis of computational studies \cite{Poole,Gallo,Palmer,Debenedetti,Gartner,Singh}. The LDL and HDL phases feature different structural and topological properties (see section \ref{HDLLDL}) and several studies in the literature proposed a variety of chemical order parameters to investigate these two forms, \cite{Poole2,Shiratani,russo2014,tanaka2019,montes2020,muthachikavil2022,Martelli,foffi2021,Shiratani2,Singh,Errington,Martelli2, montes2020,Martelli3}. We also mention works of Bak\'o et al, \cite{Bako1,Bako2}, in which the authors use a spectral clustering method, based on the analysis of the eigenvalues of the graph Laplacian, to study the global properties of an H-bonded network. 

In \cite{Benzi4}, we presented a new order parameter based on a centrality measure, the total communicability (TC), see 
\cite{Benzi}, which was found to be effective in differentiating the molecules in the LDL phase from the molecules in the HDL phase. We used for our analysis a water model that was previously shown to exhibit a metastable liquid–liquid critical point in deeply supercooled conditions and we chose four temperatures along an isobar that, in this water model, crosses the coexistence line, i.e., along which a first order phase transition between LDL and HDL occurs. 
We compared the distributions of the total communicability with the distributions of the eigenvector centrality and the degree, and we found that the last two are not well-suited to the task of characterizing molecules in the two phases. Regarding the eigenvector centrality, this can be explained by observing that the adjacency matrices associated 
with water networks exhibit a small spectral gap, hence using only the dominant eigenvector to analyze the network leads to a loss of structural information. The degree, on the other hand, being a local measure of connectivity, is unable to capture indirect  interactions between non-neighboring water molecules and thus to account for medium to long range effects.
We also computed the fraction of LDL population at  four different reference temperatures of physical interest and found a sudden and marked decrease in the LDL fraction upon raising the temperature, consistently with a transition between the two pure phases, as expected by crossing the coexistence line. 

In the present work we investigate the performance of additional centrality measures when used to identify if water is in the LDL or HDL phase and to characterize the internal organization of the two liquid forms. We consider the closeness, betweenness, Katz, and subgraph centrality, and we compare them with the results obtained with the total communicability. 
Furthermore, we also analyze the structural and connectivity properties of the two phases using other well-established metrics, such as clustering coefficients, a bipartivity measure, the algebraic connectivity, 
the graph energy, and the relative number of closed cycles of certain lengths. 

The paper is organized as follows. In section 2, we briefly describe the LDL and HDL phase; in section 3, we recall some basic concepts from  graph theory and specify the centrality measures and global metrics used in this work. In section 4 we describe the used methods of molecular dynamics simulations, in section 5, we explain how we construct the graph given a box with water molecules, while in section 6 we present the results  of our numerical experiments. Conclusions and ideas for future works are given in section 7.

\section{Low- and high-density liquid phases} 
\label{HDLLDL}

Water is a complex liquid with anomalous properties: for example, to cite a few, liquid water is denser than solid water (ice), and this implies that the ice floats instead of sinking; it has very high specific heat, meaning that the water temperature rises more slowly than the temperature of almost any other substance; liquid water at ambient pressure shows a maximum in density at around 4 °C. 
One of the most popular hypotheses for explaining many water anomalies is based on the existence of a transition between two liquid phases, referred to as low-density liquid (LDL) and high-density liquid (HDL), \cite{Poole,Gallo,Singh,Palmer,Debenedetti,Gartner,Foffi,zanetti2021}. The existence of a liquid-liquid phase transition (LLPT) between these two liquids was hypothesized in the {\em supercooled region}, 
in which the temperature is lowered below the water freezing point, without the liquid becoming a solid. The measurements in this region are challenging since supercooled liquids quickly crystallize. In contrast, glassy water at these low-temperatures is stable, and hence well characterized, and presents two structural forms: the low-density amorphous (LDA) form and the high-density amorphous (HDA) one. For this reason, it was conjectured that a counterpart of these
two forms exists in the supercooled region. 

The structural characterization of liquids is extremely difficult because, differently
from crystalline solids, liquids lack long-range order, while retaining short- to medium-range order. At a local level, the low-density liquid form of water exhibits tetrahedral order, while the high-density liquid one is more disordered and less tetrahedral, with a high probability of having interstitial water molecules within the first hydration shell. Nevertheless, their characterization on a medium-to-long length-scale is very challenging and is still a matter of debate \cite{Martelli,russo2014, Foffi}.

\section{Graph theory background and notation}\label{sec2}

Most of the material in this section is standard, but we include it for the sake of setting the notation and terminology used in the paper. 
A graph $G=(V,E)$ consists of a finite set  $V = \{v_1,...,v_N \}$ of nodes (or vertices) and a set $E \subseteq V \times V$ whose elements are called the edges (or links) of $G$. If the edges have an orientation, i.e., $(v_i, v_j)$ is a directed edge from the node $v_i$ to the node $v_j$, then the graph is {\em directed}. 
Otherwise, it is {\em undirected}, and we do not distinguish
between $(v_i,v_j)$ and $(v_j, v_i)$. A positive number (the {\em weight}) is associated with each edge. If these weights are equal to one for all the edges, then the network is called {\em unweighted}, otherwise it is {\em weighted}.

The structural information of a graph is contained in the associated adjacency matrix $A$. It is a matrix of size $N \times N$, where $N$ is the number of nodes of $G$, and the entry $a_{i, j}$ is equal to the weight of the edge $(v_i, v_j)$. If such an edge does not exist, $a_{i,j}$ is set to zero. $A$ is symmetric if and only if the associated graph is undirected. In this work we consider undirected, unweighted graphs without self-loops, i.e., without edges of the form $(v_i, v_i)$. Hence, the adjacency matrix $A$ is symmetric, binary, and with all zeros on the main diagonal.

The {\em degree }of a node $v_i$ is the number of links incident upon the node. It is equal to the sum of the entries in the $i$th row of the matrix $A$. 
A graph is {\em regular} if all the nodes in it have the same degree. 
A {\em walk} between the nodes $v_i$ and $v_j$ is a sequence of nodes that starts from the node $v_i$, follows a sequence of edges in the graph $G$, and ends at node $v_j$.
 A {\em closed walk} is a walk where the final node coincides with the initial node. A {\em path} is a walk that never visits the same vertex twice, and a {\em cycle} is a closed walk where each node is visited only once, except the initial node  (which is equal to the final one). 
A {\em shortest path} between two nodes is a path of minimal length between such nodes; the number of edges comprising a shortest path defines the 
{\em geodesic distance} between the two nodes (if no such path exists, we say that the distance between the two nodes is infinity).
 A graph $G$ is {\em connected} if for any pair $v_i, v_j \in V$, there exists a walk starting at node $v_i$ and ending at node $v_j$. In this case, the adjacency matrix $A$ is 
 {\em irreducible}, i.e., it cannot be brought to block diagonal form by a symmetric permutation of its rows and columns. A  {\em connected component} of
 a graph $G$ is a connected subgraph $G' = (V', E')$ of $G =(V, E)$ with $V' \subseteq V$ and $E'\subseteq E$.
 
Using the adjacency matrix $A$, we can obtain the number of walks in the associated graph. 
Indeed, for any positive integer $k$,  the entry $[A^k]_{i,j}$ is the number of walks of length $k$ between nodes $v_i$ and $v_j$, where $A^k$ denotes the $k$th power of $A$. 
In particular, the diagonal entries of $A^k$ count the number of closed walks of length $k$ in $G$.

Let $K$ be the {\em degree matrix} of $G$, defined as the diagonal matrix whose $i$th diagonal entry is the degree of the $i$th node.  The {\em Laplacian matrix} of $G$ is defined as
\begin{align}\label{laplacian}
L = K - A \,.
\end{align}

Note that $L$ is a singular, symmetric positive semidefinite matrix with $L \mathds{1} = 0$, where $\mathds{1}$ denotes the vector of all ones.  The multiplicity of $0$ as an eigenvalue of $L$ is equal to the number of connected components of $G$.

The {\em density} $\delta(G)$ of a graph $G$  is the ratio between the number $|E|$ of its edges and the number of possible edges, which is $\frac{N(N-1)}{2}$ for an undirected,
loopless graph with $N$ nodes. The {\em  internal density} of a subgraph  $G'=(V',E')$ of $G$ is the ratio between the number of internal edges of $V'$ and the number of all possible internal links, while the {\em external density } is the ratio between the sum external degree of $V'$, defined as the total number of edges connecting nodes in $V'$ with nodes in
$V\setminus V'$, and the number of all possible such edges. Finally,
 a {\em cluster} of nodes  is a set of connected nodes  whose internal density (as a subgraph of $G$) is larger than its external one.

\subsection{Centrality measures}

Loosely speaking, a centrality measure (CM) is a quantity that measures the importance of a node in a graph. From a mathematical point of view, we can think of
a CM as a function from the set of nodes to the set of positive real numbers, $CM : V \longrightarrow \mathbb{R}_{\ge 0}$.  In recent years, centrality measures have found 
several applications in chemistry.  For example, eigenvector centrality has been applied to the study of
certain molecular substructures in \cite{Pari}, whereas closeness and betweenness 
centrality have been used to show that osmolytes act as hubs in the  hydrogen-bond network of the solution \cite{Sundar2021}.

In this section we describe some of the most important centrality measures in use.  
We are not attempting here to present an exhaustive list of centrality measures, but we select only the most important ones from our point of view. We refer to  \cite{Estrada}
or \cite{Newman} for a more complete description.  

\subsubsection{Degree centrality} The simplest centrality measure is the degree centrality, defined simply as the degree of a node. Since it considers only the nearest neighbors, it is a purely local notion and cannot capture interactions between non-neighboring nodes.

\subsubsection{Closeness centrality} Closeness centrality (CL) measures how close a node is to the rest of the nodes in the 
graph, see \cite{Freeman79centralityin}. In the case of a connected graph, it is computed as the average inverse of the sum of the 
distances between the node and the other nodes in the graph: 
\begin{align}\label{CL}
CL(v_i) = \frac{N - 1}{s(v_i)},
\end{align}
where $s(v_i) = \sum_{v_j \in V} \{ \text{length of the shortest path between } v_i \text{ and } v_j, \,\, j\ne i \}$ and $N$ is the number of nodes in the graph $G$. 
Conventionally, if a vertex $v_i$ is an isolated node (there are no walks that reach the node), then we set $CL(v_i) = 0$.

If the graph $G$  is not connected, we consider for each node $v_i$ the connected components $G'= (V', E')$ that contain such node. Let $N' \le N$ be the number of vertices in $G'$.
Then, see \cite{wassermanfaust},
\begin{align}\label{CL_no_connected}
CL(v_i) =  \frac{(N' - 1)^2}{s(v_i)(N-1)}\,.
\end{align}
We note that if $N'= N$, we obtain equation (\ref{CL}).

\subsubsection{Betweenness centrality} Like the closeness centrality, the betweenness centrality (BC) is based on shortest paths in the network,
see \cite{Freeman1977Se}. 
In this case, the importance of a node depends on its capacity to connect different parts in the graph. Hence, a node has a high value of BC if it enables and eases the passage of information between the other nodes in the network. In the simplest version, it considers the shortest path between the vertices, and it is defined as
\begin{align}\label{BET}
BC(v_i) = \sum_{\substack{k,l = 1 \\ k \ne i \ne l}}^N \frac{ \sigma(v_k, v_i, v_l)}{\sigma(v_k, v_l)}\,,
\end{align}
where $\sigma(v_k, v_i, v_l)$ is  the number of shortest paths between $ v_k$ and $v_l$  that go through the node $ v_i $, and $\sigma(v_k, v_l)$ is the number of all shortest paths connecting $ v_k $ with $ v_l$.

Some variations of the betweenness centrality exist, where maximum flows or random walks are used instead of the shortest paths. In this work, we use only the classical version, and we refer to \cite{Estrada} for more information.

\subsubsection{Katz centrality} In many situations, in order to better understand the roles of nodes in a graph, it is advantageous to consider the effect of all walks (not just the shortest paths).  Intuitively, we want to give more weight to the shorter walks, which can be achieved by means of an  attenuation factor $\alpha > 0$. 
We recall that the powers of the associated adjacency matrix give information about the number of walks on the graph. The Katz centrality (KC)  is the first 
walk-based measure that we consider, see \cite{Katz}. It is defined as 
\begin{align}\label{Katz1}
KC(v_i) = [(I_N  + \alpha A  + \alpha^2 A^2 + \cdots + \alpha^k A^k + \cdots)\mathds{1}]_i \,,
\end{align}
with $I_N$ the $N\times N$ identity matrix and $\mathds{1} $ the vector of all ones.  In order for the matrix series in parentheses to converge, the
parameter $\alpha$ must be such that $0 < \alpha < 1/{\rho(A)}$, where $\rho(A)$ is the spectral radius of the adjacency matrix. This guarantees that
the series converges to the nonnegative matrix $(I_N - \alpha A)^{-1}$ and that  we can rewrite the expression as
\begin{align}\label{KZ}
KC(v_i) = [(I - \alpha A)^{-1}\mathds{1}]_i \,.
\end{align}

Given that here we consider undirected graphs only, the adjacency matrix $A$ is symmetric. Hence, using the Spectral Theorem we can also
express this centrality measure in terms of the eigenvectors and eigenvalues of $A$. In particular, for the Katz centrality we have
\begin{align}\label{KZ_spectral}
KC(v_i) = \sum_{k=1}^N \frac{1}{1- \alpha \lambda_k} ({\bf p}_k^T \mathds{1}){ \bf p}_k  (i) \,,
\end{align}
where $\lambda_1 \ge \lambda_2 \ge \ldots \ge \lambda_N$ are the
eigenvalues of $A$ and ${\bf p}_k = [{\bf p}_k (1), \ldots , {\bf p}_k (N)]^T$ is the eigenvector associated with $\lambda_k$, normalized with respect to the Euclidean norm.
It is easy to see (cf.\cite{Benzi2}) that for $\alpha \to 0+$, the ranking obtained by Katz centrality is the same as the one obtained using degree centrality.

\subsubsection{Eigenvector centrality} The Eigenvector Centrality (EC) is a popular centrality measure which can also be interpreted in terms of walks on the graph \cite{Bonacich}.
 %
%
 Given that we are assuming that the graph $G$ is connected, the adjacency matrix is irreducible. It is also a symmetric non-negative matrix, therefore by the Perron-Frobenius Theorem 
 \cite{HJ1}
the maximum eigenvalue $\lambda_1$ is positive and simple, and it coincides with the spectral radius $\rho(A)$. Furthermore, there exists a unique, up to normalization, eigenvector  ${\bf p}_1$ associated to $\lambda_1$, such that ${\bf p}_1 > \bf{0}$ (entry-wise). Hence, the principal eigenvector can be used to define a centrality measure. We define $ EC(v_i) = \mathbf{p}_1 (i)$, where ${\mathbf p}_1 (i)$ indicates the $i$th entry of the vector ${\mathbf p}_1$. It is not difficult to see that
\begin{align}\label{EG1}
EC(v_i) = \lim_{k \to \infty} \frac{ \# \, \text{walks of length $k$ through } v_i}{ \# \, \text{walks of length $k$ in $G$}}\,,
\end{align}
showing that a node has high eigenvector centrality if it keeps being visited ``often" by walks on $G$, as the length of the walks tends to infinity.
Moreover, it can be shown that eigenvector centrality is a limiting case of Katz centrality, in the sense that the ranking of the nodes obtained by
Katz centrality as $\alpha \to 1/\rho(A)$ tends to the ranking given by eigenvector centrality, see \cite{Benzi2}.   

\subsubsection{Subgraph centrality} The subgraph centrality (SUB) quantifies the importance of a node by considering the number of closed walks in the graph that  pass
through  such node, again giving more weight to shorter walks \cite{Estrada3}. Let $\beta > 0$, then the subgraph centrality of node $v_i$ is defined in terms of the 
matrix exponential as
\begin{align}\label{SUB}
SUB(v_i) = \sum_{k = 0}^{\infty} \frac{\beta^k}{k!}[A^k]_{i ,i} = [e^{\beta A}]_{i , i}.
\end{align}
Longer walks are penalized through the use of the rapidly decreasing weights, $\frac{\beta^k}{k!}$ for a walk of length $k$. We note that the constant $\beta$ is a tuning
parameter that can be 
used to give more or less weight to longer walks. The default value of $\beta$ is $1$. The subgraph centrality can be obtained from the spectral decomposition of the 
adjacency matrix $A$:
\begin{align} \label{SUB1}
SUB(v_i) = \sum_{k=1}^N e^{\beta \lambda_k} {\mathbf p}_k^2 (i)\,.
\end{align}
The rankings obtained with SUB reduce to the degree and eigenvector centrality rankings for $\beta \to 0+$ and $\beta \to \infty$, respectively \cite{Benzi2}.

\subsubsection{Total communicability} The Total Communicability (TC) is closely related to the subgraph centrality. Unlike SUB,
the TC does not consider only closed walks, but all the walks between a node $v_i$ and every vertex in the graph, including $v_i$ itself \cite{Benzi,Estrada4}. 
The same weighting is used as in SUB, i.e, each walk of length $k$ is  weighted by $\frac{\beta^k}{k!}$. The TC centrality of node $v_i$ is defined as
\begin{align}\label{TC}
TC(v_i) = [e^{\beta A} \mathds{1}]_i = \sum_{k = 0}^{\infty} \frac{\beta^k}{k!}[A^k \mathds{1}]_{i} \,,
\end{align}
where $\beta$ is a positive parameter. Unlike the subgraph centrality, the TC can be computed efficiently even for very large graphs using algorithms 
for evaluating the action of a matrix function on a vector \cite{BB}. Like the other walk-based measures, 
it can be expressed in terms of the eigenvalues and eigenvectors of the matrix $A$:
\begin{align}\label{TC_spectral}
TC(v_i) = \sum_{k=1}^N e^{\beta \lambda_k} (\mathbf{p}_k^T \mathds{1})\mathbf{p}_k(i)\,.
\end{align}

The limiting behavior for $\beta \to 0+$ and for $\beta \to \infty$
 of the ranking obtained using TC is the same as that of the subgraph centrality, see again \cite{Benzi2}.

\subsection{Other metrics}
A number of global measures exist that can be used to characterize the geometrical organization of a graph. In the context of water molecules, for example, it is often useful to investigate the presence of cycles of various lengths, see \cite{Foffi}.  Moreover, as is the case for other types of networks, it is of interest 
to study the clustering coefficient and the average shortest path lengths as indicators of possible small-world patterns \cite{Newman}.

\subsubsection{Counting the cycles in a graph}

Some analytic formulae based on the adjacency matrix's spectral moments exist to compute the number of cycles in a graph. We recall that the $k$th moment is 
defined as $\mu_k =\sum_{j=1}^N \lambda_j^k = \text{Tr}(A^k)$, and that $\text{Tr}(A^k)$ is the number of closed walks of length $k$ in the graph. In particular, let $|S_k|$ be the number of cycles of length k, then
\begin{align}\label{cycles}
\begin{split}
|S_3| & = \frac{1}{6} \mu_3 \,,  \\
|S_4| & = \frac{1}{8} \biggl ( \mu_4 - 2 \sum_{i=1}^N k_i (k_i-1) - 2 m \biggl)\,,  \\
|S_5| & = \frac{1}{10} \biggl ( \mu_5 - 30 |S_3| - 10 \sum_{k_i >2} \frac{[A^3]_{i ,i}}{2} (k_i-2) \biggl) \,,
\end{split}
\end{align} 
where $k_i$ is the degree of node $v_i$, $m$ is the number of edges in the graph and $[A^3]_{i, i}$ is the entry $(i,i)$ of the the third power of $A$. Note that the quantity $\frac{[A^3]_{i ,i}}{2}$ is equal to the number of triangles attached to the node $v_i$. See \cite{Alon} and \cite{Estrada} for details in obtaining the formulae.

\subsubsection{Average Watts-Strogatz clustering coefficient and transitivity index}

The clustering coefficient and the transitivity index give information about how clustered a network is locally and globally, respectively
\cite{Estrada,Newman}. The clustering coefficient of node $v_i$, denoted $C_i$, is defined as 
\[
C_i  = \frac{2 t_i}{k_i (k_i -1)}\,,
\]
where $t_i$ is the number of triangles that pass through the node $v_i$ of degree $k_i$. It represents the fraction of actual triangles over all potential triangles 
passing through node $v_i$. 
The average Watts-Strogatz clustering coefficient is just the mean of the clustering coefficients over all nodes:
\[
\bar{C} = \frac{1}{N} \sum_{i=1}^N C_i\,.
\]
The transitivity index (or Newman clustering) is a global measure of the frequency of triangles in a network. It is defined as the number of triangles in the graph divided by the total number of walks of length 2, $| P_2 |$:
\[
C = \frac{3 | S_3 |}{| P_2 |} \,.
\]

Note that $|P_2| = \sum_{i=1}^N k_i(k_i - 1)/2$.
Both these global indexes take values between 0 and 1; values close to 0 indicate a low density of triangles in the graph, while values close to 1 represent a 
high density. Generally speaking, the two clustering measures are in a good correlation, but there are graphs where they behave differently, see \cite{Estrada5,Benzi3}.

\subsubsection{Bipartivity measure}

A graph is bipartite if we can partition the node-set into two disjoint non-empty sets $V_1$ and $V_2$, such that every edge of the graph connects a vertex in $V_1$ to one in $V_2$. This means that in a bipartite graph, the closed walks of odd length are not present. A bipartivity measure quantifies how much a graph is close to being bipartite. Estrada has
introduced the following measure:
\begin{align} \label{Bip}
B = \frac{\text {Tr}(\cosh(A))}{\text {Tr}(\exp(A))}\,,
\end{align}
where $\cosh (A) = \frac{\exp(A )+ \exp(-A)}{2}$ denotes the hyperbolic cosine of the adjacency matrix $A$  \cite{Estrada2}.
We note that $0.5 < B \le 1$, and that $B=1$  if and only if the graph is bipartite. Values of $B$ close to 1 indicate that the graph
is nearly bipartite (i.e, it can be made into a bipartite graph by rewiring a small number of edges).

\subsubsection{Average shortest path length}

Recall that the distance $d(v_i, v_j)$ between the nodes $v_i$ and $v_j$ is the number of edges in a shortest path between such vertices. If a walk between $v_i$ and $v_j$ does not exist, their distance is set to $+ \infty$. The diameter of a graph is defined as the maximum distance between any two vertices in the graph:
\begin{align}\label{diameter}
\text{diam}\, (G) = \max_{v_i, v_j \in V} \{ d(v_i, v_j) \}\,.
\end{align}

We can define the distance matrix $D$ as the symmetric matrix whose  $( i, j )$ entry  is equal to $d(v_i, v_j)$, for all $i,j = 1,\ldots ,N$. Then the average shortest path length
(ASPL) is defined as
\begin{align}\label{ASPL}
\bar l = \frac{\mathds{1}^T D \mathds{1}}{N (N-1)}\,,
\end{align}
where we consider only pairs of nodes for which a path connecting them exists (such a path always exists if $G$ is connected).  A relatively large clustering coefficient and an
average path length that grows slowly (say, logarithmically) with the number of nodes $N$ characterize the class of small world (or Watts--Strogatz) networks \cite{Newman}. In contrast,
very small clustering coefficients and a fast growing  (with $N$) ASPL are observed for grid-like graphs and for classical random graph models. 

\subsubsection{Energy of a graph}

Let $\lambda_1, \lambda_2, \ldots , \lambda_N$ be the eigenvalues of the adjacency matrix $A$, then the energy of the associated graph $G$ is defined as
\[
E(G) = \sum_{i=1}^N |\lambda_i| \,.
\]
The energy is an important graph invariant inspired by chemistry (H\"uckel Molecular Orbital Theory) and that has been extensively studied; see, e.g, the
monograph \cite{Gutman}, as well as
\cite{EB} for an interpretation of the graph energy in terms of walks on $G$ and for applications to molecular structure.  It is interesting to note that adding edges to 
an existing graph does not always lead to an increase in $E(G)$; see, e.g., \cite{Arizmendi}.

\subsubsection{Algebraic connectivity}

Let us denote the eigenvalues of the Laplacian matrix of $G$ as follows:
\begin{align}\label{AlCo}
0 = \mu_1 \le \mu_2 \le \ldots  \le \mu_n \,.
\end{align}
The second eigenvalue of the Laplacian matrix, $\alpha(G) = \mu_2$, is called the {\em algebraic connectivity} of the
graph $G$, and provides a quantitative measure of ``how well-connected” the graph is \cite{Fiedler}; note that $\alpha (G) > 0$ if and only if $G$ is connected. Larger values of $\alpha (G)$ correspond to graphs that are not easily disconnected by removing a small number of edges. It is known \cite{Fiedler} that the algebraic connectivity is edge-monotone, i.e., for two graphs with the same vertex set, $G_1 = (V,E_1)$ and $G_2 = (V,E_2)$, with $E_1 \subseteq E_2$, it is $\alpha(G_1) \le \alpha(G_2)$.

\subsubsection{Degree assortativity coefficient} 
This quantity measures the tendency of a node to connect with nodes that have equal, or almost equal, degrees \cite{Newman}. If this quantity is positive, then the network is assortative. Otherwise, it is called disassortative. It is defined as
\begin{align}\label{assortativity}
r = \frac{|P_2| (|P_{3/2}| +C - |P_{2/1}|)}{3|S_{1,3}| +|P_2| (1 - |P_{2/1}|)} \,,
\end{align}
where $|P_i| $ is the number of paths of length $i$, $|P_{i/j}|:=|P_i| /|P_j| $, $|S_{1,3}|$ is the number of star fragments of four nodes and $C$ is the transitivity index
defined earlier. The denominator is always grater than or equal to zero, therefore $r > 0$ (the graph is assortative) if, and only if, $|P_{2/1}| < |P_{3/2}| +C$. The denominator is equal to zero in the case of a regular graph; in this case, the assortativity coefficient is not defined.

\section{Molecular dynamics simulations of liquid water}\label{sec3}
Classical molecular dynamics (MD) is a powerful method to sample the configurational space of a liquid (or, more generally, of a solution). Here, we use previously performed \cite{zanetti2021} MD simulations of pure water; in particular, we use MD simulations on the nanosecond timescale of TIP4P/2005 water at 1950 bar and at four different temperatures, namely 170 K (200 ns), 180 K (150 ns), 200 K (100 ns) and 240 K (100 ns). 
The TIP4P/2005 water model \cite{abascal2005general}, which is widely used to simulate  both pure water and aqueous solutions \cite{zanetti2021,Zanetti19,Biswas21}, exhibits a metastable liquid–liquid critical point in deeply supercooled conditions.
Given the simulation parameters, which are consistent with those of Biddle et al. \cite{biddle2017two}, the liquid-liquid critical point is estimated to be at 1700 bar and 182 K. A rectangular box containing 710 water molecules was used and the MD simulations were performed in the NPT ensemble with the 5.1.2 version of the GROMACS software \cite{abraham2015gromacs}.
Temperature and pressure were kept constant by using the velocity rescaling temperature coupling \cite{bussi2007canonical} and the Parrinello–Rahman barostat with 2 ps relaxation times \cite{parrinello1980crystal}. Periodic boundary conditions were used, long range electrostatic interactions were treated with the particle mesh Ewald method \cite{darden1993particle} with a real space cutoff of 0.9 nm and for short range interactions a cut-off radius of 0.9 nm was employed.  
All bonds and intramolecular distances were constrained using the LINCS algorithm \cite{hess1997lincs} and a 2 fs time step was utilized.

\section{Liquid water as a network}

Given a MD trajectory, the coordinates of a given set of atoms extracted at a given time frame can be used to build a network. Since in the simulations periodic boundary conditions are imposed, meaning that the box is replicated along the three directions of space, each molecule interacts with the images of the molecules which are on the opposite sides of the box.

We form the graph $G = (V, E)$ where the water molecules represent the nodes of $G$ while links between the molecules are the edges. There are different ways of defining a link between two particles, \cite{Bondar}. In this work, using similar ideas and arguments to our previous paper \cite{Benzi4}, we consider two molecules as connected by an edge if the distance between their oxygen atoms is less than or equal to $0.35$ nm (considering their replications along the three directions). In order to evaluate the distance, we compute the matrix of the physical distances $D_{phy}$. In practice, for each pair of particles $v_i$ and $v_j$ in the box,  we calculate their distance along each direction. For example, let $dist_x(v_i,v_j)$ be their distance along the x-axis. If this quantity is larger than half of the box size in that direction $L_x$, we replace it with the length of the box minus this value: $dist_x(v_i,v_j) = L_x - dist_x(v_i,v_j)$. We do the same for each direction, and we determine the physical distance matrix:
\[ 
D_{phy} (v_i,v_j) = \sqrt{dist_x(v_i,v_j)^2 + dist_y(v_i,v_j)^2 + dist_z(v_i,v_j)^2} \,.
\]
At this point, if $D_{phy} (v_i,v_j) \le 0.35$, there exists an edge between the nodes $v_i$ and $v_j$. The corresponding adjacency matrix $A$ is symmetric, undirected, unweighted, and of dimension $N \times N$, where $N$ is the number of molecules (in our experiments $N = 710$). 

The code we used to run all network analyses, together with instructions on how to install and use it, has been made available in the GitHub repository (https://github.com/ChiaraFaccio/WaterNetworks). The code makes heavy use of the NetworkX 2.6.3 module \cite{nx} in Python 3.7, the Python package  NetworkSNS \cite{bucci2},  the code \cite{pbc} for imposing the periodic boundary conditions, and the open-source packages NumPy \cite{harris2020array}, SciPy \cite{2020SciPy-NMeth}, and Matplotlib \cite{Hunter:2007}. All the computations were performed on a laptop with a 4 Intel core  i7-8565U CPU @ 1.80GHz  - 1.99 GHz and 16.0GB RAM.

\section{Results}\label{sec4}
We analyze here four MD simulations on the nanosecond timescale of neat water using the TIP4P/2005 water model \cite{Singh} along the 1950 bar isobar \cite{zanetti2021,Benzi4}, which is above the critical pressure (1700 bar) and along which the liquid–liquid coexistence line is crossed at a temperature of around 175 K. We chose temperatures below (170K) and above (180, 200 and 240 K) this liquid-liquid transition temperature in order to maximize changes in water structure on going from the lowest temperature, at which the water is in the LDL-form, to the highest ones, at which the HDL-form prevails.
For each MD trajectory, we extract 100 frames equally-spaced in time and analyze the corresponding networks using the previously described graph- and
network-theoretic metrics. We note explicitly that while the number of nodes  in each network remains constant, the number of edges varies along the MD trajectory as well as with the 
temperature, since the spatial arrangement of the oxygens varies.

\begin{figure}[!h]
\centering
\subfloat[]{\includegraphics[width=2.5in]{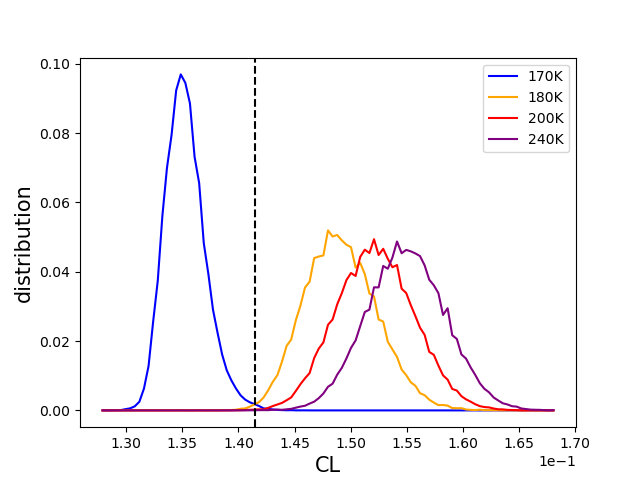} \label{a}}
\hfil
\subfloat[]{\includegraphics[width=2.5in]{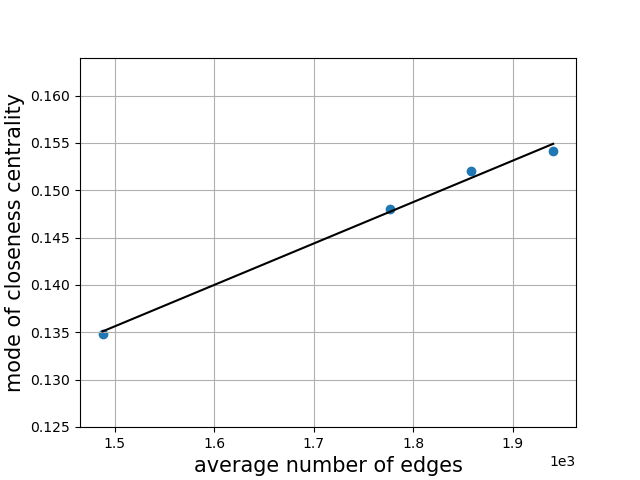}
\label{b} }
\vfil
\subfloat[]{\includegraphics[width=2.5in]{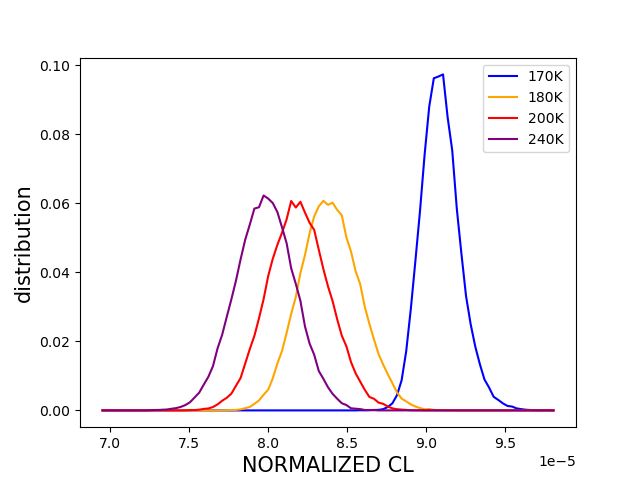} \label{c} }
\hfil
\subfloat[]{\includegraphics[width=2.5in]{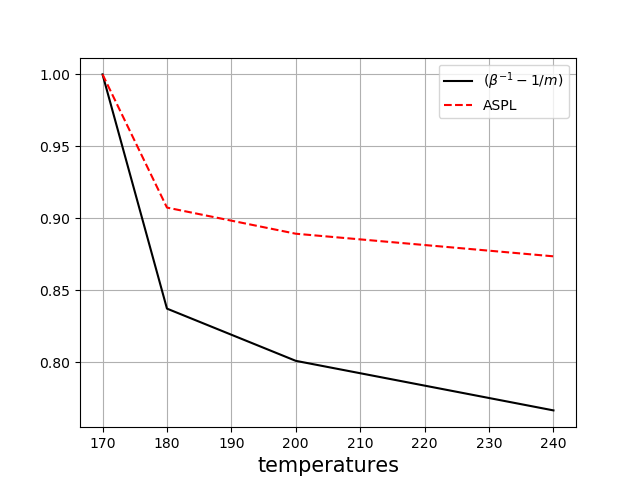}
\label{d} }
\caption[]{Figure[\ref{a}]: distribution of closeness centrality values. The black dashed line in the left graph corresponds to the intersections between the distributions at 170K and the distributions at 180K; Figure[\ref{b}]: linear correlation between the average number of edges in the four temperatures and the mode of CL for the four distributions; Figure[\ref{c}]: distribution of the closeness centrality values normalized by the number of edges; Figure[\ref{d}]: average shortest path length and average $(\beta(G)^{-1} - 1/m )$ 
coefficient for the four temperatures. The data are obtained from the MD simulations at 1950 bar at 170K (blue), 180K (orange), 200K (red) and 240K (violet).}
\label{FIG1}
\end{figure}

\subsection{Analysis of water networks via centrality measures}

We begin by comparing the distributions for the four temperatures using different centrality measures to see how well they can distinguish the LDL phase from the HDL phase. 
Degree and eigenvector centrality, as mentioned in the Introduction (see also \cite{Benzi4}), are not effective for this task and therefore we will not show the corresponding results, while the distributions of total communicability are reported in Figure[\ref{atc}].
In Figure[\ref{a}], we report the distributions of closeness centrality at 1950 bar. This centrality measure is able to differentiate the LDL phase from the HDL one: the distribution at 170K is well-separated from the distributions at higher temperatures, while there is a significant overlap between the plots at 180K, 200K, and 240K. Closeness centrality values follow the same trend as the number of links in the network. In Figure[\ref{b}], we can see that there is, in fact, a linear correlation between the modes of the distributions at the four temperatures and the average number of edges in the graphs. When the temperature increases, the distributions are shifted towards higher values of CL. Since we construct the adjacency matrices using a threshold on the physical distances between the oxygen atoms, the distributions indicate that in the HDL phase, the average distance between the water molecules is smaller than in the LDL phase.

Some comments on the effect of normalization of the computed centrality values are in order, 
since we are comparing networks with the same number of nodes but different number of links. 
If we normalize by the number of edges, we obtain the distributions in Figure[\ref{c}]. 
The normalized closeness values now show an opposite behavior, i.e., when the temperature increases, the modes of the distributions are shifted towards smaller values. This 
apparently surprising
behavior can be attributed to the difference between the rates of change of the number of edges and of the quantity $s(v_i)$ in equation
(\ref{CL})  for the four temperatures. Let $m$ be the number of edges in the network, then we have
\begin{align}\label{CL_norm_edges}
\frac{CL( v_i)}{m} = \frac{N-1}{m}\frac{1}{s (v_i)} = \biggl (\frac{N}{m} - \frac{1}{m} \biggl )\frac{1}{s (v_i)} = \biggl (\beta(G)^{-1} - \frac{1}{m} \biggl )\frac{1}{s (v_i)}\,,
\end{align}

where $\beta(G)= m/N$ is the {\em beta index} of the graph $G$ (not to be confused with the parameter $\beta$ used in the definition of sugraph and total
communicability centrality), an indicator of the level of connectivity of a graph. In our case, the average quantity $(\beta(G)^{-1} - 1/m )$, reporting on the increase of the number of edges upon raising the temperature, is equal to 0.48, 0.40, 0.38, and 0.37, on going from 170 K to 240 K. As a global indicator for the behavior of $s(v_i)$, we can consider the average shortest path length (ASPL). We obtain that ASPL is 7.38, 6.70, 6.57, and 6.45 at 170 K, 180 K, 200 K and 240 K, respectively. In Figure[\ref{d}] we plot these values, normalized to be between 0 and 1, and we observe that the rate of decrease of $(\beta(G)^{-1} - 1/m )$ is greater than the rate of the ASPL. For this reason, normalizing by the number of edges we obtain the distributions of Figure[\ref{c}].

\begin{figure}[!h]
\centering
\subfloat[]{\includegraphics[width=2.5in]{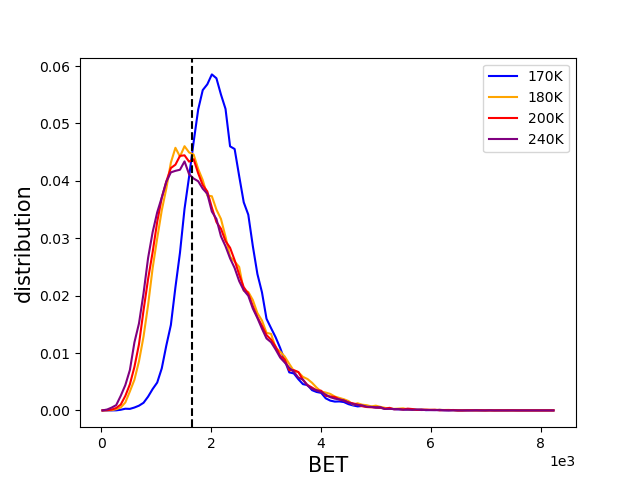} \label{abet}}
\hfil
\subfloat[]{\includegraphics[width=2.5in]{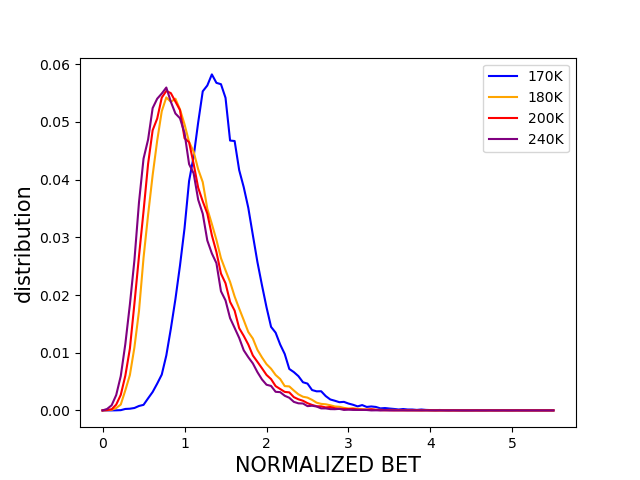}
\label{cbet} }
\caption[]{Figure[\ref{abet}]: distribution of betweenness centrality values. The black dashed line in the left graph corresponds to the intersections between the distribution at 170K and the distribution at 180K;  Figure[\ref{cbet}]: distribution of the betweenness centrality normalized by the number of edges. The data are obtained from the MD simulations at 1950 bar at 170K (blue), 180K (orange), 200K (red) and 240K (violet).} 
\label{FIG2}
\end{figure}

The distributions of the betweenness centrality values are reported in Figure[\ref{abet}]. 
We observe that the BC values in the LDL phase tend to be somewhat higher than for the HDL phase. This can be explained by noting
that the BC of a node is proportional to the number of shortest paths that pass through such vertices. Since in the high-density liquid water the molecules are closer, and thus the number of edges increases, there is a rise in the availability of shortest paths, and the fraction of shortest paths through each vertex tends to decrease. However, there is a significant overlap between the two regimes, and betweenness centrality does not appear to be able to clearly identify the two density forms of liquid water.
The normalization by the number of links (Figure[\ref{cbet}]) determines only minor changes (i.e, the distributions at higher temperatures become tall like at 170K), and the BC is still unable to well separate between the LDL and HDL phases. We conclude that the betweenness centrality is not an useful  order parameter to identify the two density forms of liquid water.

\begin{figure}[!h]
\centering
\subfloat[]{\includegraphics[width=2.5in]{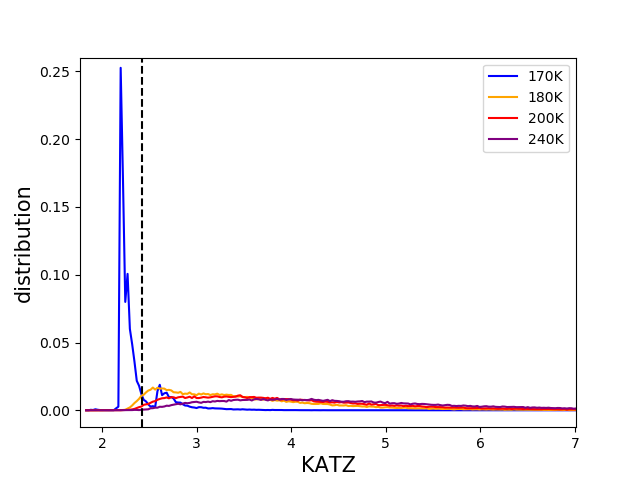} \label{akatz}}
\hfil
\subfloat[]{\includegraphics[width=2.5in]{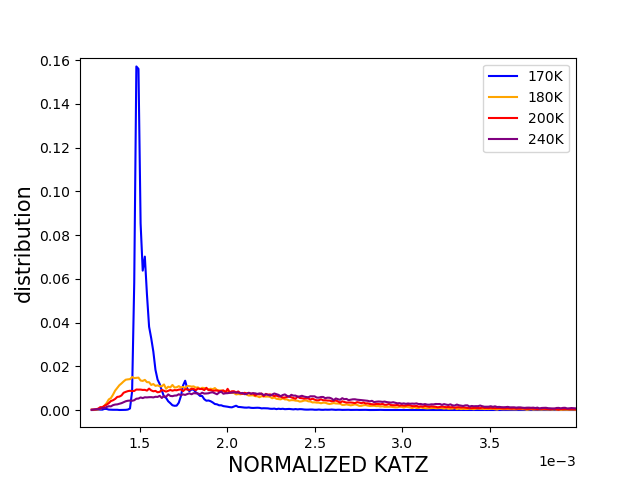}
\label{ckatz} }
\caption[]{Figure[\ref{akatz}]: distribution of Katz centrality values using $\alpha = 0.1362$. The black dashed line in the left graph corresponds to the intersections between the distribution at 170K and the distribution at 180K;  Figure[\ref{ckatz}]: distribution of the Katz centrality normalized by the number of edges. The data are obtained from the MD simulations at 1950 bar at 170K (blue), 180K (orange), 200K (red) and 240K (violet).}
\label{FIG3}
\end{figure}

\begin{figure}[!h]
\centering
\subfloat[]{\includegraphics[width=2.5in]{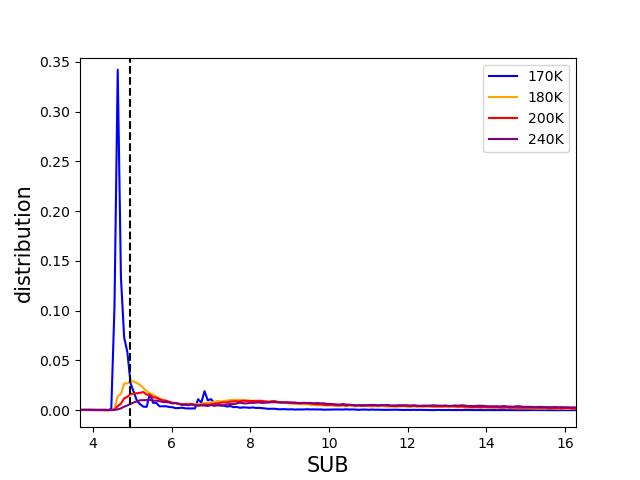} \label{asub}}
\hfil
\subfloat[]{\includegraphics[width=2.5in]{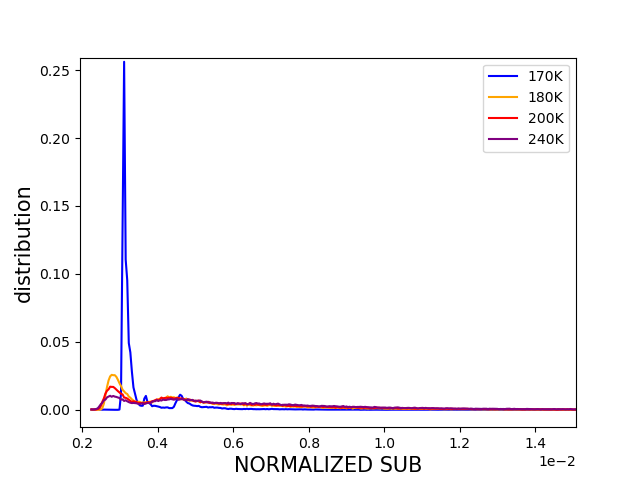} \label{csub} }
\caption[]{Figure[\ref{asub}]: distribution of subgraph centrality values using $\beta = 1$. The black dashed line in the left graph corresponds to the intersections between the distribution at 170K and the distribution at 180K;  Figure[\ref{csub}]: distribution of the subgraph centrality normalized by the number of edges. The data are obtained from the MD simulations at 1950 bar at 170K (blue), 180K (orange), 200K (red) and 240K (violet).}
\label{FIG4}
\end{figure}

\begin{figure}[!h]
\centering
\subfloat[]{\includegraphics[width=2.5in]{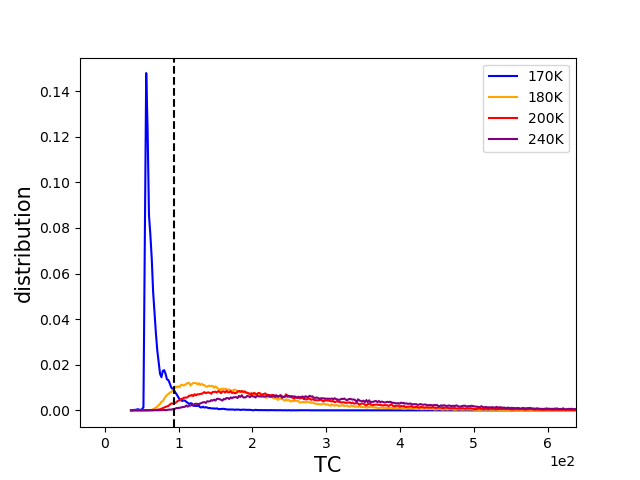} \label{atc}}
\hfil
\subfloat[]{\includegraphics[width=2.5in]{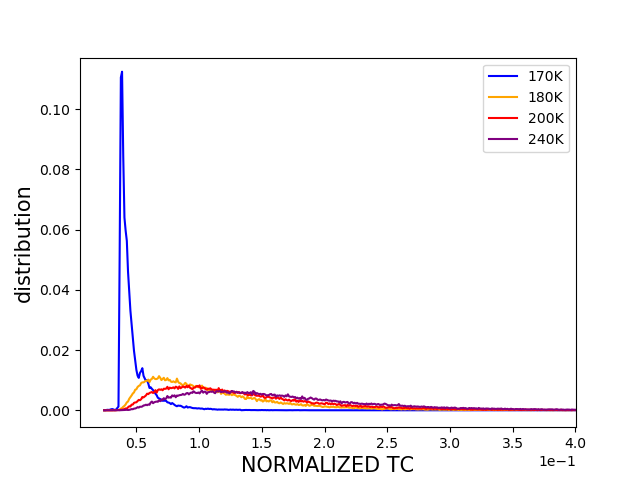} \label{ctc} }
\caption[]{Figure[\ref{atc}]: distribution of total communicability values using $\beta = 1$. The black dashed line in the left graph corresponds to the intersections between the distribution at 170K and the distribution at 180K;  Figure[\ref{ctc}]: distribution of the total communicability normalized by the number of edges. The data are obtained from the MD simulations at 1950 bar at 170K (blue), 180K (orange), 200K (red) and 240K (violet).}
\label{FIG5}
\end{figure}

\begin{figure}[!h]
\centering
\subfloat[]{\includegraphics[width=2.6in]{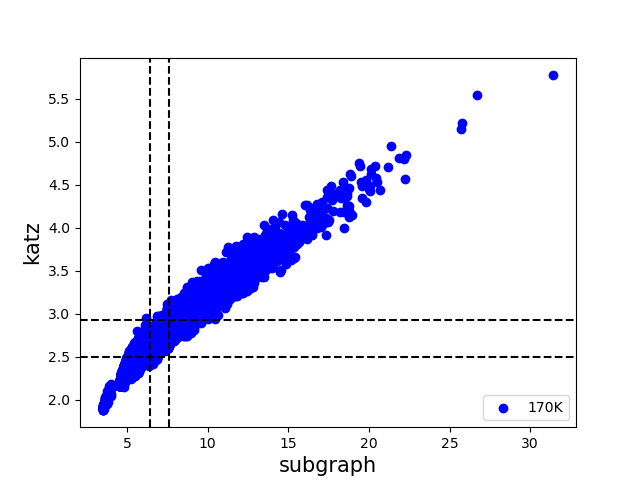} \label{subkatz}}
\hfil
\subfloat[]{\includegraphics[width=2.6in]{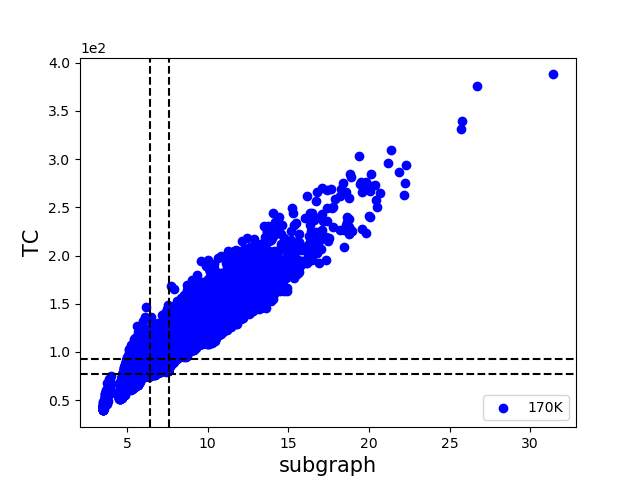}\label{subtc}}
\vfil
\subfloat[]{\includegraphics[width=2.6in]{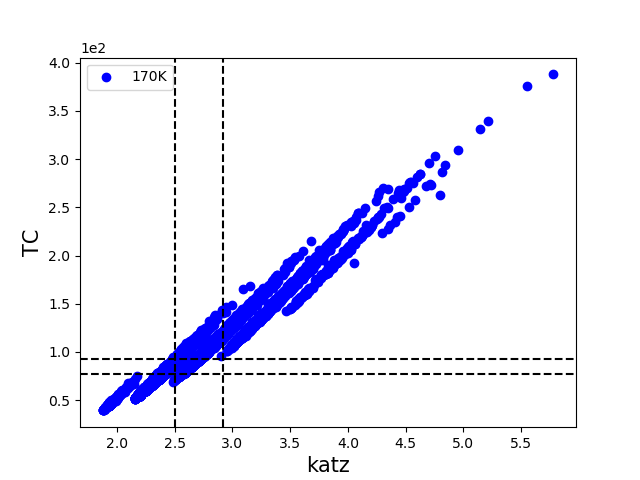}\label{tckatz}}
\caption[]{Figure[\ref{subkatz}]: correlation between the values of SUB and KC for the temperature at 170K; Figure[\ref{subtc}]: correlation between the values of SUB and TC for the temperature at 170K; Figure[\ref{tckatz}]: correlation between the values of TC and KC for the temperature at 170K. The black dashed lines correspond to the bounds for the peaks.}
\label{FIG6}
\end{figure}

The distributions of the Katz, subgraph, and total communicability centrality values are reported in Figures [\ref{akatz},\ref{asub},\ref{atc}], respectively. Recall that
the KC of a node, like the TC, considers all the walks through that node, while SUB counts only the closed walks. All these centralities use a parameter to give less weight to the long walks. These measures have similar distributions, and they well identify the LDL phase from the HDL one. In particular, for Katz centrality we chose the penalty factor $\alpha = 0.1362$, which we determined as follows. Recall that $\alpha$ is subject to the condition
$0 < \alpha < 1/ \rho(A)$, where $\rho(A)$ is the spectral radius of the adjacency matrix, and that if $\alpha \to 1/\rho(A)$, the node ranking of KC converges to the ranking of EC, while, for $\alpha \to 0$, the ranking reduces to the ranking of the degree centrality. The spectral radii in the LDL phase are lower than in the HDL phase, so if we compute the KC of each graph using a penalty factor that depends on the spectral radius of the associated adjacency matrix $A_i$, $\alpha = 1/(\gamma \cdot \rho(A_i))$ (with $\gamma > 1$ a constant), then, as the temperature increases, we give less weight to the walks. On the other hand, in the HDL phase, there are more walks of length $k$ than in the LDL phase for $k \ge 1$, but since $\alpha$ is smaller, we tend to lose this information. For this reason, we set the $\alpha$ parameter to be the same for all four trajectories. The maximum spectral radius for all the graphs is 6.67. As a consequence,  the $\alpha$ factor must be at least less than $1/6.67 \approx 0.15$, and far enough from this quantity to avoid reducing to the eigenvector centrality. We obtain that $\alpha = 1/ (1.1 \cdot 6.67) = 0.1362$ is a good choice, in fact the centrality measure allows to identify the two densities of the liquid water. Nevertheless, at 170 K (LDL phase) there is a small secondary peak overlapping with the distributions obtained at the higher temperatures (i.e., in the HDL phase). This peak is also present in the subgraph centrality distribution at 170K with the choice $\beta = 1$. The molecules belonging to the two peaks are almost the same, as we can observe from Figure[\ref{subkatz}], and they correspond, in part, to the molecules of the small peak in the TC distribution at 170K with TC values between 77 and 93.

For the subgraph and total communicability centrality, we employ $\beta =1$ for two reasons: first of all, it is widely used as a default value and, secondly, following the results in \cite{Benzi4} for the total communicability, this choice leads to including the contributions of molecules within a 1 nm distance (the contribution of molecules at a greater distance
being numerically negligible), which is the typical correlation length in the pair radial distribution function for the $O \cdot \cdot \cdot O$ contacts, see \cite{Benzi4}.

If we normalize these centralities by the number of edges, KC is no longer able to identify the two density phases, subgraph centrality distributions show an opposite behavior (similar to closeness),  while for TC, the overlap between the distributions increases slightly. In particular, TC is less influenced by this normalization than Katz and Subgraph. This is probably due to the fact that the three normalized centrality measures give different weights to the walks: for example, at 170K with an average of $m=1487.61$ edges, given $\alpha = 0.1362$ and $\beta = 1$, the penalty factor referred to the degree is $9.16 \cdot 10^{-05}$, $3.36 \cdot 10^{-04}$, $6.72 \cdot 10^{-04}$ for Katz, Subgraph and Total Communicability respectively.
In practice, normalizing by the number of edges leads to a loss of information about an essential feature of the two density forms, namely, the fact that
water molecules in the HDL phase display a higher connectivity. 

Next, we calculate the LDL fraction at the four temperatures using the previous values without normalization. As in \cite{Benzi4}, 
we use the crossing point between the distributions at 170K and 180K at 1950 bar to discriminate the two phases. Let $x_*$ be such point, then $v_i$ is selected in the LDL phase if $CM(v_i) \le x_*$, otherwise it is assigned to the HDL one.

\begin{table}[!h]
\caption[]{Fraction of LDL population at 1950 bar for the four temperatures as obtained using different centrality measures.}
\label{TABLE1}
\centering
\begin{tabular}{|c|c|c|c|c|}
\hline
 & 170K & 180K & 200K & 240K \\
\hline
CL & 99.53 \% & 0.27 \% & 0.03 \% & 0.009 \% \\ 
BC & 84.18 \% & 61.22 \% & 59.16 \% & 57.18 \%\\
KC & 79.81 \% & 3.15 \% & 0.79 \% & 0.18 \%\\
SUB & 71.27 \% & 8.51 \% & 3.80 \% & 1.49 \% \\
TC & 89.43 \% & 6.48 \% & 1.74 \%& 0.20 \% \\
\hline
\end{tabular}	
\end{table}

In Table[\ref{TABLE1}], the results are summarized. It can be observed that, with the exception of the betweeness, all centrality measures well discriminate the two liquid phases. In the case of the betweeness, as previously recognized in Figure[\ref{abet}], there is a significant overlap between the distributions at the four temperatures, and thus BC overestimates the fraction of LDL water molecules in the pure HDL regime.

With respect to CL and TC, the Katz and subgraph centrality provide a lower estimate of the LDL fraction at 170 K. This is due to the presence of the already commented secondary peaks. In the case of the total communicability, the secondary peak is located before the crossing point used to define the two regimes and does not overlap with the distributions obtained in the HDL phase. 

A strong point of centrality measures is that they provide information on each node of the graph. Therefore, as described in \cite{Benzi4}, they work at 
the molecular level and permit to obtain information about the internal organization of the network. 

To this aim, we focus on Katz, subgraph, and total communicability centrality
(we exclude the betweenness centrality because it does not identify the two liquids and the closeness because the nodes with high centrality values are all in the center of the box. Its definition introduces in fact a bias towards central nodes and therefore it cannot be used to investigate the structural organization of the water network).
 When the temperature increases, the distributions of Katz, subgraph and TC become broader, meaning that the values of the centrality measures tend to be higher. We previously observed \cite{Benzi4} that in the HDL phase, the nodes with high TC are not isolated, but appear to organize themselves into patches. 

\begin{figure}[!h]
\centering
\subfloat[Katz centrality]{\includegraphics[trim=4.5cm 1cm 4cm 1cm, clip=true, width=2.8in]{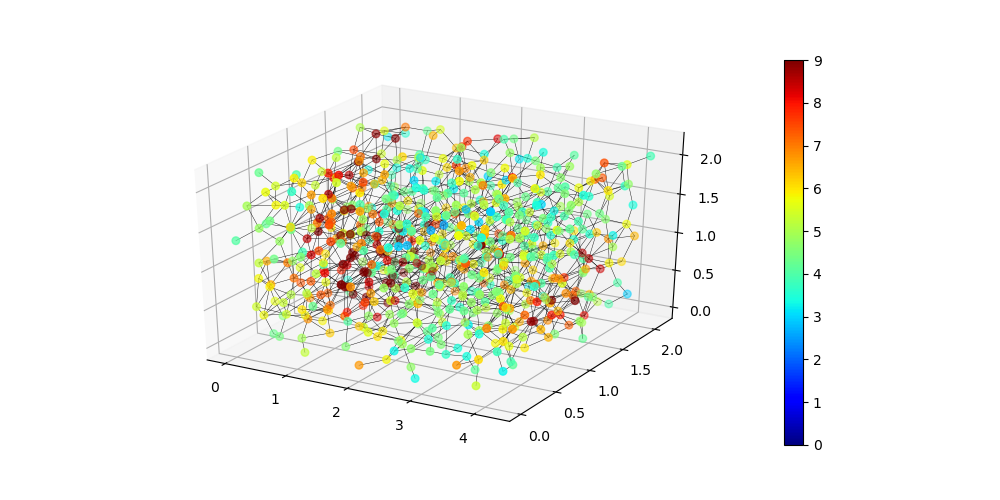} \label{dist_katz}}
\hfil
\subfloat[Subgraph centrality]{\includegraphics[trim=4.5cm 1cm 4cm 1cm, clip=true, width=2.8in]{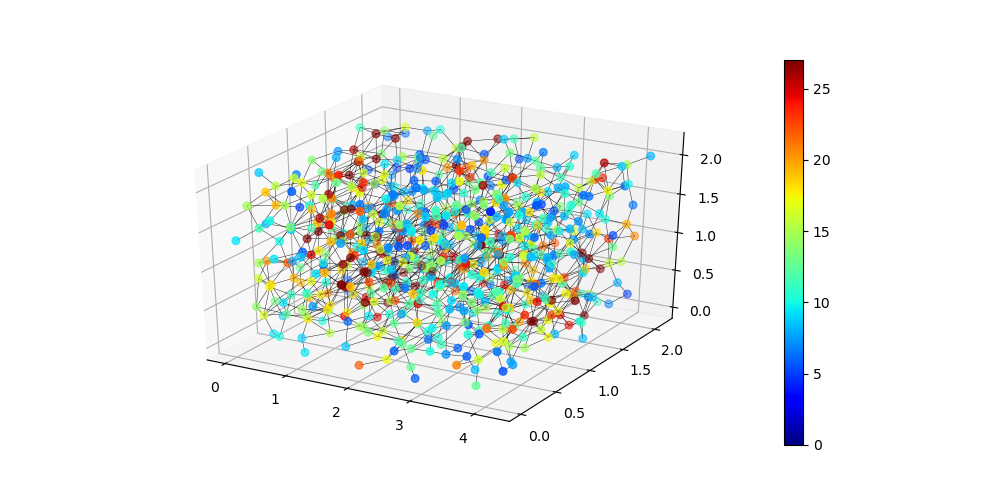} \label{dist_sub}}
\vfil
\subfloat[Total communicability]{\includegraphics[trim=4.5cm 1cm 4cm 1cm, clip=true, width=2.8in]{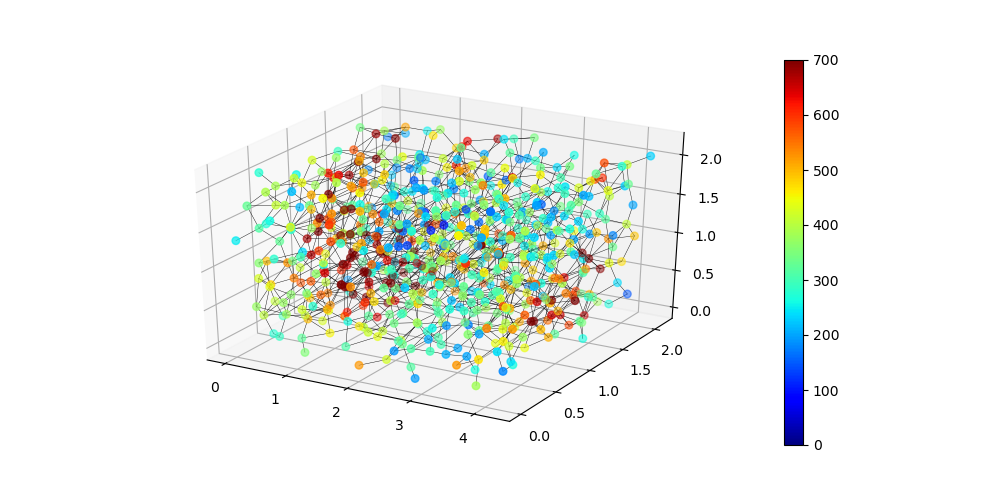} \label{dist_TC}}
\caption{3D  plot of the distribution of oxygen atoms with different Katz, subgraph and TC centrality values for frame 7 at 240K.}
\label{FIG7}
\end{figure}

We consider the trajectory at 240K, where almost all the water molecules are in the HDL phase. From Figure[\ref{dist_katz}] and Figure[\ref{dist_sub}], we observe the same behavior that we noted in \cite{Benzi4} also for Katz and subgraph centrality, i.e., molecules with high centrality measures are assembled in specific regions of the box.

\begin{figure}[!h]
\centering
\subfloat[]{\includegraphics[trim=2.5cm 1cm 2.5cm 1cm, clip=true, width=2.8in]{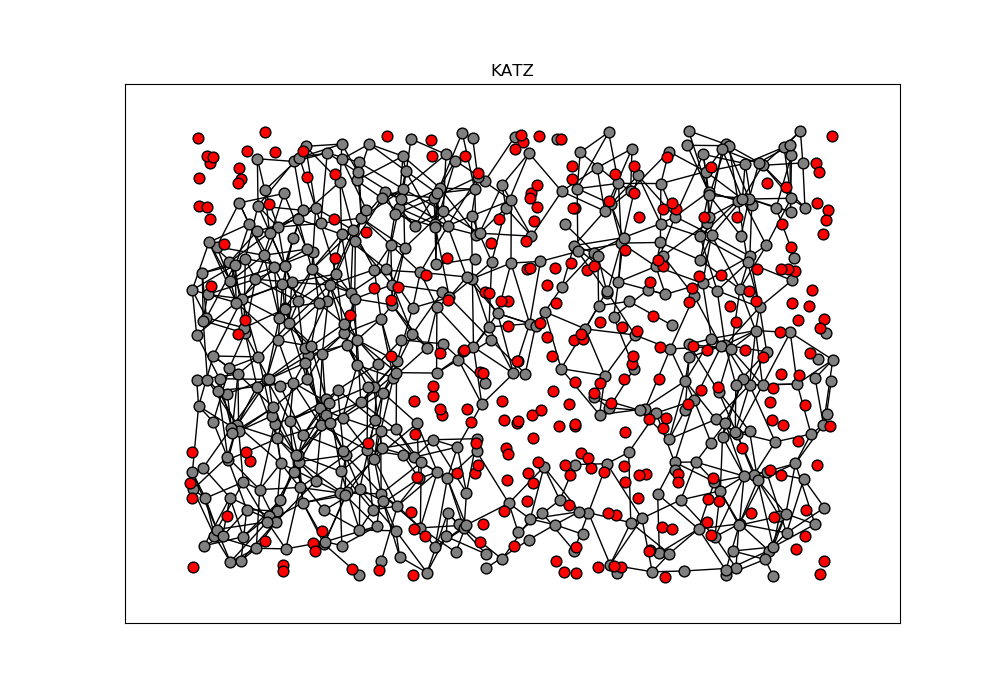} \label{katz7}}
\hfil
\subfloat[]{\includegraphics[trim=2.5cm 1cm 2.5cm 1cm, clip=true, width=2.8in]{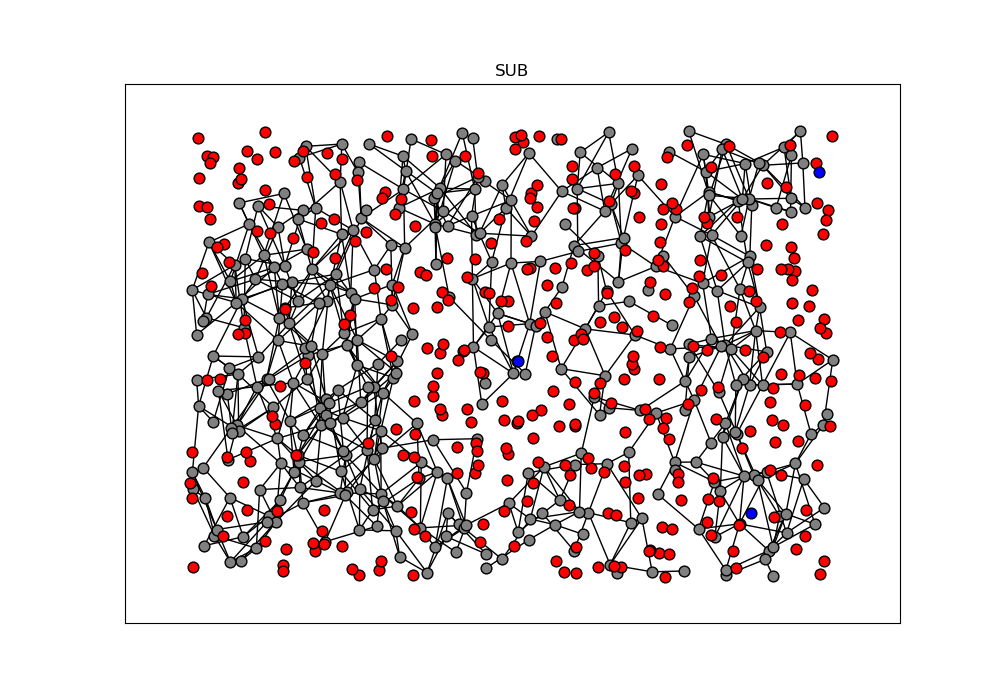} \label{sub7}}
\vfil
\subfloat[]{\includegraphics[trim=2.5cm 1cm 2.5cm 1cm, clip=true, width=2.8in]{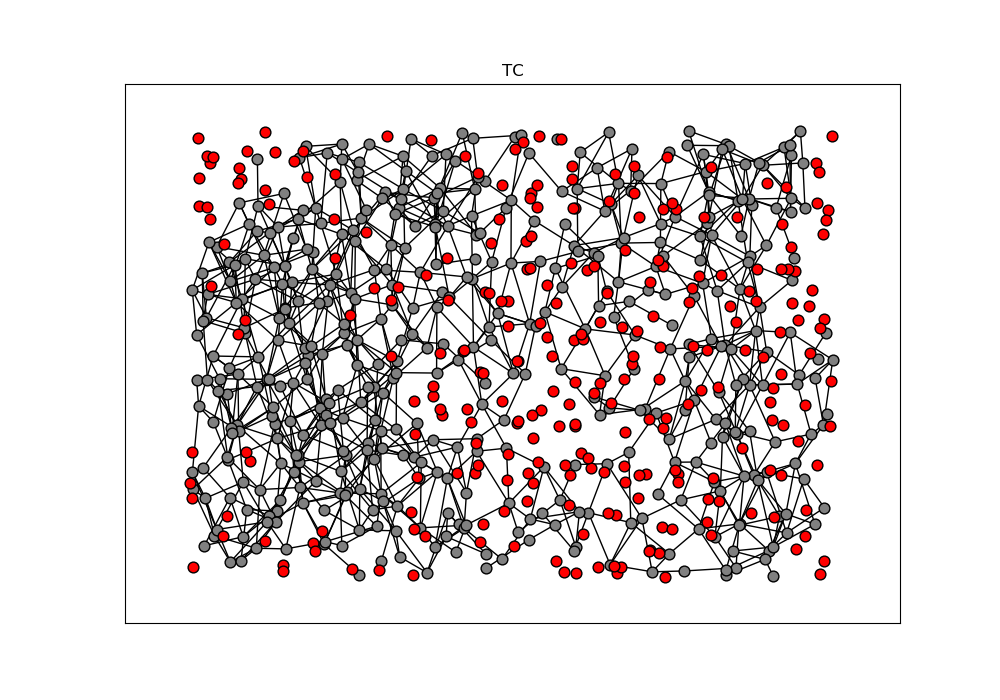} \label{tc7}}
\caption[]{Snapshots of the internal organization of the water at 240K, frame 7. Blue nodes represents the LDL molecules, red  nodes the HDL molecules and the silver ones are the molecules with CM above the average. The edges highlight the connections among high-CM nodes.}
\label{FIG8}
\end{figure}

As shown in Figure[\ref{FIG8}], there is a big patch of molecules with centrality values above the average (colored in silver in the Figure) for all three centrality measures. The red nodes are the HDL molecules, while the blue ones are LDL molecules. Only the subgraph centrality identifies some of the molecules as being in the low-density phase, according to results in Table[\ref{TABLE1}]. The regions are composed of 447 molecules in the TC case, 448 molecules in the Katz case, and 360 for SUB. We argue that the patch size obtained with SUB differs from those found with KC and TC because SUB does not consider all the walks of all the lengths (like KC and TC), but only the closed walks. The average internal density of the patches obtained with the three centrality measures is 0.01, the average external one is $7.29 \cdot 10^{-5}$, while the density  $\delta(G)$ of the entire graph is $8.05 \cdot 10^{-3}$. We underline that the above values refer to the graph density, which is not necessarily correlated to the physical density of the water network. We also mention that we tried to apply standard
community detection techniques to our networks, such as the Girvan--Newman algorithm or the greedy modularity technique (see \cite{Fortunato} for details), but we obtained different clusters which, however, do not appear to be amenable to chemical-physical interpretation. We plan to further investigate this aspect in future works. 

\subsection{Computational aspects}
In a typical simulation, a large number of water networks need to be analyzed, and computational efficiency is paramount. 
In Table[\ref{TABLE2}], we report the average execution times for computing the centrality of every node for a single frame (network). 
For the closeness and betweenness centrality, the computation of all the shortest paths is by far the more expensive part. It requires a breadth-first search (BFS), which needs $O(|V| + |E| )$ time for each node of the graph. The Katz centrality can be evaluated using a direct solver (sparse Cholesky factorization) for the equation (\ref{KZ}). 
The subgraph centrality looks at the diagonal entries of the exponential of the adjacency matrix, and it can be computed using quadrature rules and the Lanczos algorithm
(see \cite{BB} for details). On the other hand, the total communicability corresponds to the product of the exponential of the adjacency matrix times the vector of all ones, and it can be computed efficiently using a Lanczos-based algorithm for evaluating the action of a matrix function on a vector \cite{BB}.

\begin{table}[!h]
\centering
\caption{Average of the execution times of a frame at different temperatures for the graph with 710 nodes using different centrality measures. The values are in seconds. }	
\label{TABLE2}
\begin{tabular}{|c|c|c|c|c|c|}
\hline
& CL & BC & KC & SUB & TC \\
\hline
170K & 0.72 & 2.64 & 0.0063 & 13.64 & 0.0062 \\
180K & 0.78 & 2.84 & 0.0082 & 16.79 & 0.0054 \\
200K & 0.79 & 2.64 & 0.0085 & 16.74 & 0.0054 \\
240K & 0.81 & 2.72 & 0.0080 & 16.10 & 0.0060 \\
\hline
\end{tabular}
\end{table}

From the table, we can see that the TC centrality is the fastest measure to compute, followed by the Katz centrality, while SUB is the most time-consuming. 

\subsection{Results using global metrics}

Here we analyze the two liquid phases using some global metrics from  graph and network  theory. 
Water networks are 3D grid-like networks which in the LDL phase are near to being regular graphs (with common degree 4).  As the temperature increases,
this regular structure is perturbed as some degree of disorder sets in.  However, even in the HDL phase the networks remain fairly regular, and do not display the emergence
of the small world phenomenon. This is easily explained by noting that, because of the threshold we use to define an edge, the appearance of new links is limited to pairs of molecules that are not too far in space,
hence the connectivity tends to remain short-ranged. 
 In Table[\ref{TABLE3}], we display the average shortest path length (ASPL) and diameter 
 for the four different temperatures. We observe that both ASPL and diameter exhibit a moderate reduction for increasing temperature. 
In the third column of the table we report the (average) value of the algebraic connectivity of the networks for the four different temperatures. As expected, the algebraic connectivity increases monotonically with the temperature.

Next, we compare our water networks with the classical Erd\"{o}s--R\'enyi random graph model, in which each pair of nodes is connected with a given probability $p>0$.
In an Erd\"{o}s--R\'enyi random graph, the density of $G$ is equal to its clustering coefficient, see \cite{Estrada}. Since for both the  LDL and HDL phases the densities $\delta(G)$ are low while the average Watts--Strogatz indexes (see Table[\ref{TABLE4}])
are relatively high, their clustering is not similar to that of an Erd\"{o}s-R\'enyi random graph, reflecting the much higher transitivity of water networks.
Note that as the temperature increases the graph density also increases, reflecting the creation of new links between pair of molecules as the network becomes less ordered.

\begin{table}[!h]
\caption{Mean values of the average shortest path length, diameter, algebraic connectivity, density, bipartivity measure, and energy  of water network for the four temperatures. }	
\label{TABLE3}
\centering
\begin{tabular}{|c|c|c|c|c|c|c|}
\hline
 & ASPL & $\text{ diam} (G)$ &  $\alpha (G)$ & $\delta(G)$ &  bipart. & $E(G)$ \\
\hline
170K & 7.38 & 13.20 & 0.101 & 0.0059 & 0.94 & 1271.14 \\ 
180K & 6.70 & 12.04 & 0.132 & 0.0071 & 0.80 & 1356.80 \\
200K & 6.57 & 12.00 & 0.138 & 0.0074 & 0.77 & 1379.63 \\
240K & 6.45 & 11.94 & 0.150 & 0.0077 & 0.73 & 1402.80 \\
\hline
\end{tabular}
\end{table}

In ice, oxygen atoms form a hexagonal crystalline lattice, where each water molecule is bonded to four other water molecules, in a tetrahedral arrangement, through hydrogen bonds. 
(For a 4-regular graph, the number of edges is $m =N\cdot 4/2$; if $N = 710$, then $m = 1420$. The average number of links at 170K is 1487.61, very close to this value).
Such a regular lattice is easily shown to be bipartite. Since liquid water at 170K behaves similarly to ice, we expect that for water networks at this low temperature
the bipartivity measure defined in (\ref{Bip}) is close to 1 (nearly bipartite graph).  At higher temperatures, when this ordered structure is perturbed by the appearance of new links (see first column of Table[\ref{TABLE6}]),
the (near) bipartivity of the network is lost, and we expect the bipartivity measure to decrease. This intuition is confirmed by the results reported in Table[\ref{TABLE3}]. This behaviour is confirmed by inspecting the adjacency spectrum at different temperatures.  Let $\lambda_1 \ge \lambda_2 \ge \ldots \ge \lambda_N$ be the eigenvalues of $A$. It is a basic fact of spectral graph theory that a graph is bipartite if and only if   $\lambda_N = -\rho(A) = -\lambda_1$. In Table[{\ref{TABLE5}}], we report the average extreme values of the spectrum of $A$ and the mean values of their sum. The negative eigenvalue changes only by a small amount, while the larger one grows as the temperature increases. As a result, their sum moves away from 0. This behavior is also reflected in the increase of the graph energy $E(G)$ with the temperature, see Table[\ref{TABLE3}]. Hence, we have another confirmation that the graph becomes less bipartite as the temperature increases.

\begin{table}[!h]
\caption{Mean values of the minimum $\lambda_N$ and maximum $\lambda_1$ eigenvalue of the adjacency matrices $A$, and the mean values of their sum  $\lambda_1 + \lambda_N$ for the four temperatures.}	
\label{TABLE5}
\centering
\begin{tabular}{|c|c|c|c| }
\hline
 & $\lambda_N$ & $\lambda_1$ & $ \lambda_1 + \lambda_N$  \\
\hline
170 K & -3.54 & 4.58 & 1.04 \\
180 K & -3.83 & 5.74 & 1.90 \\
200 K & -3.88 & 5.93 & 2.05 \\
240 K & -3.91 & 6.13 & 2.22 \\
\hline
\end{tabular}
\end{table}

Next, we have found that the Watts-Strogatz clustering coefficients $\bar{C}$ are in a linear correlation with the transitivity indexes $C$
(see Table[\ref{TABLE4}]), with a Pearson correlation of 0.97 for all the temperatures considered. All the networks are assortative due to the fact that $|P_{2/1}| < |P_{3/2}| +C $. In particular, the quantity  $|P_{2/1}|$ is larger than the average of $|P_{3/2}|$, $0 \le C \le 1$ hence the transitivity index $C$ has a primary effect on the assortativity. As the temperature increases, the gap $ (|P_{2/1}| - |P_{3/2}|)$ grows, therefore  the transitivity index also increases. The increase in the transitivity index upon raising the temperature, related to the increase in the number of edges in the graph, corresponds, from a chemico-physical point of view, to the increase in the number of interstitial water molecules in the high density regime (see section \ref{HDLLDL}).

\begin{table}[!h]
\caption[]{Mean values of the $|P_{2/1}|$, $|P_{3/2}|$, transitivity index $C$, Watts-Strogatz clustering coefficient $\bar{C}$, and degree assortativity coefficient $r$ for the four temperatures. }	
\label{TABLE4}
\centering
\begin{tabular}{|c|c|c|c|c|c|}
\hline
 & $|P_{2/1}|$ &$|P_{3/2}|$ & $C$ &  $\bar{C}$  & $r$ \\
\hline
170K & 3.24 & 3.23 & 0.041 & 0.029 & 0.308 \\ 
180K & 4.21 & 4.16 & 0.132 & 0.102 & 0.296 \\
200K & 4.45 & 4.38 & 0.151 & 0.122 & 0.276\\
240K & 4.70 & 4.60 & 0.175 & 0.149 & 0.252 \\
\hline
\end{tabular}
\end{table}

\begin{figure}[!h]
\centering
\includegraphics[width=2.5in]{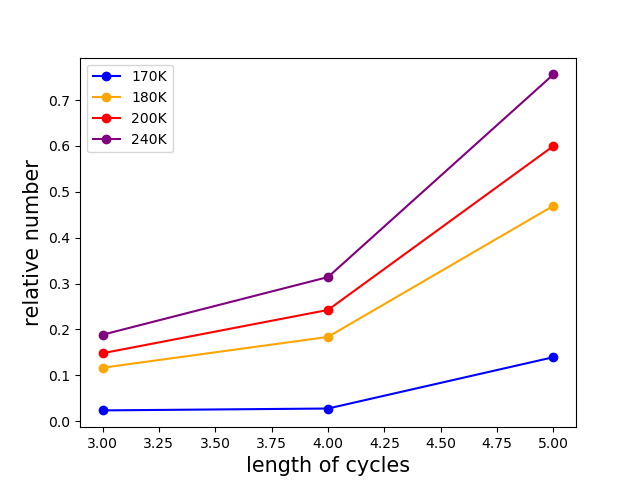}
\caption{Relative number of cycles of length 3, 4, 5 at 170K (blue), 180K (orange), 200K (red) and 240K (violet). }
\label{FIG10}
\end{figure}

Lastly, in the literature there are several works on the problem of counting the different number of cycles in the LDL and HDL phases, e.g.,  \cite{Martelli,Foffi}.
These authors use tools derived from statistics. In contrast, in this paper we use the formulas (\ref{cycles}), known from graph theory. In Figure[{\ref{FIG10}}], we observe the relative number of closed paths of length three, four, and five for the four temperatures. We normalize the number of each cycle by dividing for the maximum value among all these numbers. Finally, we report the mean at each length. In Table[\ref{TABLE6}], we also report the mean values of the number of cycles of length 3, 4, and 5.
Our results show that the number of such cycles tends to increase with the temperature, signaling a topological difference between the two phases. 

\begin{table}[!h]
\caption{Mean values of the number of edges and mean values of the number of cycles of length 3,4, and 5.}	
\label{TABLE6}
\centering
\begin{tabular}{|c|c|c|c|c| }
\hline
 & edges & cycles of length 3 & cycles of length 4 & cycles of length 5 \\
\hline
170 K & 1487.61 & 66.43 & 77.98 & 394.36 \\ 
180 K & 1776.95 & 329.36 & 519.31 & 1327.49 \\
200 K & 1857.66 & 418.64 & 685.50 & 1695.79 \\
240 K & 1940.87 & 532.69 & 888.11 & 2136.90\\
\hline
\end{tabular}
\end{table}

\section{Conclusions and future work}\label{conclusion}

In this paper, we extended the study in \cite{Benzi4} to give a structural analysis of liquid water using molecular dynamics simulations coupled with concepts and tools from graph
and network theory. We have assessed the ability of several popular centrality measures to differentiate the high-density liquid phase from the low-density liquid phase along the 1950 isobar. Since the centrality measures are computed at a molecular level, they also allow us to identify, within a single phase, the presence of high (or low) density patches.
Our comparison has confirmed the usefulness of the total communicability in this context, and found that this centrality measure is also  computationally very efficient.

Even if these networks are only moderately complex, network science tools prove extremely useful for obtaining fine information about the structure of liquid water and lead to a topological characterization of the two forms. Specifically, graph-theory based measures are able to show that the low-density phase is similar to a lattice, while in the high-density form this regularity is partly lost. These measures can also increase our physical understanding of the LDL and HDL forms. In fact, the changes we observe in both centrality measures and global metrics highlight topological differences within the HDL form. We plan to deepen this aspect in future work.

Future efforts will also be aimed at improving the description of liquid water by means of a directed graph model  in which the connections represent the hydrogen bonds among water molecules. Another improvement will be to consider a box with water molecules as a weighted graph, where the weight of an edge measures the length of the hydrogen bond. 
Finally, we would like to gain a better understanding of the clusters found by standard  techniques for community detection, which do not seem to bear much relation
with the low- and high-density patches identified by means of centrality measures such as TC.


\bibliographystyle{comnet}
%


\end{document}